\documentclass[12pt,a4paper]{article}

\usepackage{amsmath,amssymb}

\def\sq{\hbox{$\blacksquare$}}
\def\qed{\ifmmode\sq\else{\unskip\nobreak\hfil
\penalty50\hskip1em\null\nobreak\hfil\sq
\parfillskip=0pt\finalhyphendemerits=0\endgraf}\fi\vskip6mm}

\newcommand{\NNs}{\ensuremath{\mathbb{N}^*}}
\newcommand{\RR}{\ensuremath{\mathbb{R}}}
\newcommand{\RRp}{\ensuremath{{\mathbb{R}_+}}}
\newcommand{\cB}{\ensuremath{\mathcal{B}}}
\newcommand{\cF}{\ensuremath{\mathcal{F}}}
\newcommand{\cN}{\ensuremath{\mathcal{N}}}
\newcommand{\cC}{\ensuremath{\mathcal{C}}}
\newcommand{\sub}{\subseteq}
\newcommand{\s}{{^{-1}}}
\newcommand{\ag}{\alpha}
\newcommand{\bg}{\beta}
\newcommand{\cg}{\gamma}
\newcommand{\eg}{\epsilon}
\newcommand{\lm}{\lambda}
\newcommand{\LM}{\Lambda}
\newcommand{\diam}{\mathrm{diam}}
\newcommand{\proof}{\noindent\textbf{Proof. }\par\noindent}
\newcommand{\ds}{\displaystyle}
\newcommand{\und}{\underline}

\begin{document}

\parskip=3pt

\title{On Completion of Metric Mappings}
\author{Giorgio NORDO}
\date{}
\maketitle

\begin{abstract}
An internal characterization of complete metric mappings (by means
of Cauchy nets tied at a point) is given and a construction of the completion of a
metric mapping is presented.
\end{abstract}

\par\noindent\textbf{Introduction}

The notion of \textit{complete trivially metric mapping} as generalization
to mappings of the corresponding notion for metric spaces
(see, for example, \cite{fgt} and \cite{pasy}) was
recently introduced by Pasynkov in \cite{metric}.
\\
Let $X$ be a set, $Y$ be a
topological space and $f:X\to Y$ be a mapping.
\\
A pseudometric $d$ on $X$ is called a \textit{metric} (or a trivial metric
in \cite{metric}) on the
mapping $f$ if for every $y \in Y$, $d|_{f\s(\{y\}\times f\s(\{y\}}$
is a metric on $f\s(\{y\})$.
\\
In such a case, the pair $\big( f:X\to Y,\, d \big)$ is called a \textit{metric mapping}.
\\
The family $\tau(d)\cup f\s(\tau(Y))$ (where $\tau(d)$ is the topology generated
by the pseudometric $d$ on $X$ and $\tau(Y)$ is the topology of
$Y$) is a subbase for a topology on $X$ that will be called the
\textit{metric topology of the mapping $f$} and that will be denoted with
$\tau(f,d)$.
Obviously, a base for this topology is given by $\cB(f,d)=\{ U\cap f\s(V) : \, U \in \tau(d),
\, V\in \tau (Y)\}$ (see also the definition of \textit{base of a
mapping} in \cite{pasy}) and in case $X$ is a space whose topology
$\tau (X)$ coincides with $\tau(f,d)$, the mapping $f:X \to Y$ is
said \textit{trivally-metrizable} (see \cite{metric} and \cite{trivial}).
\\
A metric $d$ on $f:X\to Y$ is said \textit{complete}
if for any $y \in Y$ and any filter $\cF$ on
$X$ such that $f\s(\cN_y)\sub \cF$ and ${\forall\eg>0} \quad \exists F \in \cF \, : \,\,
\diam(F)<\eg$ then $f\s(\{y\}) \cap a_X(\cF) \ne \emptyset$
with respect to the topology $\tau (f,d)$.

In \cite{metric}, it is proved that every metric mapping $f:X\to Y$ can be
densely isometrically embedded in a complete metric mapping
$f^*:X^* \to Y$ (called the \textit{completion} of $f$)
which is unique up to isometric equivalences.
That result was obtained by reducing the proof to an equivalent
case for spaces (by using the notion of \textit{mapping parallel
to a space}).

In this paper, we will give a characterization of complete
mappings by means of the notion of Cauchy nets tied at a point.
This furnish an analogous for mappings of a well-known result
for complete metric spaces.
\\
In the second section, by using the previous result,
we will give an alternative (internal) construction
for the completion of a metric mapping.
\vskip 4mm

\par\noindent\textbf{1. Characterization of complete mappings}
\vskip 2mm

\noindent \textbf{Definition 1.}
Let $d$ be a metric on the mapping $f:X\to Y$ and $y$ a point of
the space $Y$.
A net $\Phi:D\to X$ on the set $X$ is called a \textit{Cauchy net (on the mapping $f$) tied}
to $y$ if:
\begin{itemize}
\item for any $\eg >0$ there exists some $\ag \in D$ such that for every
$\bg,\cg \ge \ag$, $d\left( \Phi(\bg), \Phi(\cg) \right)<\eg$, and
\item for any neighborhood $O$ of $y$ in $Y$ there exists
some $\ag \in D$ such that $\Phi(D_\ag) \sub f\s(O)$
(where $D_\ag = \{ \bg \in D : \, \ag \le \bg \}$ is the terminal
set of $D$ at $\ag$).
\end{itemize}

\noindent \textbf{LEMMA 2.}
Let $\big(f:X \to Y, \, d \big)$ be a matric mapping and
$\Phi:D\to X$ is a Cauchy net tied to a point $y \in Y$.
Then
$$a_X(\Phi)\cap f\s(\{ y \}) = c_X(\Phi) \cap f\s(\{ y \})$$
with respect to the metric topology $\tau(f,d)$ of the mapping $f$.

\proof
Obviously, we have only to prove that
$a_X(\Phi)\cap f\s(\{ y \}) \sub c_X(\Phi) \cap f\s(\{ y \})$.
Let $x \in a_X(\Phi)\cap f\s(\{ y \})$.
For every neighborhood $N$ of $x$ in $(X,\tau(f,d))$ there are
$U \in \tau (d)$ and $V \in \tau (Y)$ such that
$x \in U \cap f\s(V) \sub N$.
So, there is some $\eg >0$ such that $B(x,\eg)\sub U$.
Since $\Phi$ is a Cauchy net, there exists $\cg \in D$ such that
for every $\ag,bg >cg$, $d(\Phi(\ag),\Phi(\bg))<\frac{\eg}{2}$.
Since $x$ is a cluster point of $\Phi : d \to X$ with respect to
$(X,\tau(f,d))$ there exists some $\mu' > \cg$ such that
$\Phi(\mu')\in B\big(x, \frac{\eg}{2} \big) \cap f\s(V)$
and so $d(\Phi(\mu'), x) <\frac{\eg}{2}$.
Since $\Phi$ is tied to $y$, there exists some $\mu'' \in D$ such that
for every $ag \ge \mu''$, $\Phi(\ag) \in f\s(V)$.
So, there exists some $\mu \in D$ such that $\mu \ge \mu'$ and $\mu \ge \mu''$
and for every $\ag > \mu$ we have that $\Phi(\ag) \in f\s(V)$ and
\begin{equation}
\begin{split}
d\big( \Phi(\ag), x \big) & \le
d\big( \Phi(\ag), \Phi(\mu')\big) + d\big( \Phi(\mu'), x \big)
\\
& < \frac{\eg}{2} + \frac{\eg}{2} = \eg
\end{split}
\end{equation}
and so that $\Phi(\ag) \in B(x,\eg)\cap f\s(V) \sub U \cap f\s(V) \sub N$.
This proves that $x$ is a limit point of $\Phi$ with respect to
$\tau(f,d)$ and concludes our proof.
\qed

\noindent\textbf{THEOREM 3.} A metric $d$ on a mapping $f:X\to
Y$ is complete if and only if for any $y \in Y$, every
Cauchy net tied to $y$ converges to some point of the fibre
$f\s(\{y\})$ with respect to the metric topology $\tau (f,d)$ of $f$.

\proof
$(\Longrightarrow )$
Let $(f,d)$ be complete, $y \in Y$ and $\Phi:D\to X$ a Cauchy net
tied to $y$.
It is well-known (see \cite{eng}, for example) that
$$\cF (\Phi) = \{ \Phi(D_\ag) \, : \,\, \ag \in D\}$$
is a filter on $X$ that has exactly the same cluster and limit
points as $\Phi$.
\\
For every $O \in \cN_y$, since $\Phi$ is tied to $y$,
we have that there is some $\ag \in D$ such that $\Phi(D_\ag) \sub
f\s(O)$. This implies that $f\s(O)\in \cF(\Phi)$ and so that
$f\s(\cN_y)\sub \cF(\Phi)$.
Furthermore, for every $\eg >0$, since $\Phi:D\to X$ is a Cauchy
net, there exists some $\ag \in D$ such that for any $\bg,\cg
>\ag$, $d\big(\Phi(\bg), \Phi(\cg)\big)<\eg$.
Hence, $F=\Phi(D_\ag) \in \cF(\Phi)$ is such that $\diam(F)<\eg$
and by Definition of complete metric mapping, it follows that $f\s(\{y\})\cap
a_X\big( \cF(\Phi)\big)\ne\emptyset$.
Thus, $f\s(\{y\})\cap a_X( \Phi)\ne\emptyset$ and as $\Phi$ is a
Cauchy net tied to $y$, by Lemma 2, we also have that
$f\s(\{y\})\cap c_X( \Phi))\ne\emptyset$.

\noindent $(\Longleftarrow )$
Let $y \in Y$ and $\cF$ be a filter on $X$ such that
$f\s(\cN_y)\sub \cF$ and $\forall \eg >0$, $\exists F \in \cF$
such that $\diam(F)<\eg$.
It is trivial to verify that the set
$$D=\big\{ (x,F)\in X\times \cF \, : \,\, x\in F\big\}$$
is directed by the partial order $\le$ defined, for any
$(x,F),(y,G) \in D$, by setting
$$(x,F)\le(y,G) \quad \Longleftrightarrow \quad G\sub F$$
and it is a well-known fact (see \cite{eng}) that the map
$$\Phi(\cF) : D \to X \quad \mbox{ such that } \quad
(x,F) \mapsto \Phi(\cF)(x,F)=x$$
is a net on X having exactly the same cluster and limit points
as $\cF$.
\\
Now, for every $\eg >0$, there exists some $F \in \cF$ such that
$\diam(F)<\eg$.
So, there is some $x\in F$ and $(x,F)\in D$.
Hence, for every $(y,G),(z,H)\in D$ such that $(y,G),\,(z,H)\ge
(x,F)$ we have that
$d\big( \Phi (\cF )(y,G),\,  \Phi(\cF )(z,H) \big)
= d(y,z) < \eg$ because $y \in G \sub F$ and $z \in H \sub F$.
\\
Furthermore, for every $O \in \cN_y$, since $f\s(\cN_y)\sub \cF$,
it follows that $f\s(O)\in \cF$ and there is some $x \in f\s(O)$.
So, $(x,f\s(O))\in D$ and, for any $(y,G)\in D$ such that $(y,G)\ge
(x,f\s(O))$, we have that $\Phi(\cF)(y,G)=y \in G \sub
f\s(O)$ and so that
$\Phi(\cF)\big(D_{(x,f\s(O))}\big) \sub f\s(O)$.
This proves that $\Phi(\cF)$ is a Cauchy net tied to $y$.
\\
Hence, by hypothesis, we have that
$$f\s(\{y\}) \cap c_X\big(\Phi(\cF)\big)\ne \emptyset$$
and, by Lemma 2
$$f\s(\{y\}) \cap a_X\big(\Phi(\cF)\big)\ne \emptyset.$$
Since, $a_X(\cF)=a_X\big(\Phi(\cF)\big)$, we conclude that:
$$f\s(\{y\}) \cap a_X(\cF)\ne \emptyset$$
i.e. that the $d$ is a complete metric on the mapping $f:X\to Y$.
\qed
\vskip 4mm

\noindent \textbf{2. Construction of the completion of a metric mapping}

Let $d$ and $d'$ be metrics
on the mappings $f:X\to Y$ and $f':X'\to Y$ respectively.
A morphism $\lm:f \to f'$ (i.e. a continuous mapping from $X$ to $X^*$ such that
$f^*\circ \lm=f$) is called an \textit{isometric embedding} if it is
injective and $d'(\lm(x_1),\lm(x_2))=d(x_1,x_2)$ for every $x_1,
x_2 \in X$.

\noindent \textbf{Definition 4.} \cite{metric} We say that a mapping $f^*:X^* \to Y$
is a \textit{completion} of a metric mapping $f:X\to Y$ if it is a complete
metric mapping and some dense isometric embedding $\lm:f \to f^*$
is fixed.

In \cite{metric}, Pasynkov proved that every metric mapping $f:X\to Y$ has a
completion which is unique up to isometric equivalences.
The existence of the completion was proved indirectly (by means of
...) and for this reason it would be interesting to furnish a direct
construction here.
\vskip 2mm

\par\noindent
Let $d$ be a metric on a mapping $f:X\to Y$.
For any fixed $y \in Y$, let us consider the set $\cC_y$ of
all the Cauchy nets on $X$ tied to $y$ and the equivalence relation $\sim_y$ that,
for any $\Phi : D \to X$ and $\Phi : D' \to X$ in $\cC_y$, is defined by setting:
\begin{equation*}
\begin{split} \Phi \sim_y {\Phi'} \quad \Longleftrightarrow \quad
& \forall \eg>0 \, , \quad \exists \cg \in D , \,\, \exists \cg ' \in D' \, :\\
& \quad \forall \ag \ge \cg , \, \forall \ag' \ge \cg' \,\,
\Rightarrow \, d\big(\Phi(\ag), {\Phi'}(\ag')\big)<\eg
\end{split}
\end{equation*}

Let $\cC_y^*={\cC_y}_{\diagup{\sim_y}}$ be the set of all the equivalence
classes $\Phi^*$ of $\cC_y$ with respect to $\sim_y$ and let us
consider the disjoint union
$$X^*=\bigcup_{y \in Y} (\cC_y^* \times \{y\}).$$

\par\noindent
For every $(\Phi^*,y), \,( {\Phi'}^* ,y') \in X^*$, we consider the mapping
$d^*:X^*\times X^* \to \RRp$ defined by
$$d^*\big((\Phi^*,y),\, ({\Phi'}^*,y')\big)=\lim_{(\ag,\ag') \in D\times D'}d(\Phi(\ag),
{\Phi'}(\ag'))$$
where, obviously the partial order on the directed set $D\times D'$
is given by $(\ag,\ag') \le (\bg,\bg') \,\, \Longleftrightarrow \,\,
\ag\le \bg \mbox{ and } \ag'\le\bg'$.
\\
The above definition is well-posed since it that limit is unique
because $\Phi$ and ${\Phi'}$ are Cauchy nets and so
$\{ d(\Phi(\ag), {\Phi'}(\ag'))\}_{(\ag,\ag')\in D\times D'}$
is a Cauchy net on $\RR$.
\\
Now, for any $(\Phi^*, y), \, ({\Phi'}^* ,y'), \, ( {\Phi''}^* ,y'') \in X^*$,
we have that
\begin{equation*}
\begin{split}
d^*\big((\Phi^*, y) & ,\, ({\Phi''}^*,y'')\big)
 = \lim_{(\ag,\ag'') \in D\times D''}d(\Phi(\ag),
{\Phi''}(\ag''))
\\
& \le \lim_{(\ag,\ag'') \in D\times D''} \big( d(\Phi(\ag),
{\Phi'}(\ag'))+ d({\Phi'}(\ag'), {\Phi''}(\ag'')) \big)
\\
& = \lim_{(\ag,\ag') \in D\times D'} d(\Phi(\ag),{\Phi'}(\ag'))
+ \lim_{(\ag',\ag'') \in D'\times D''} d({\Phi'}(\ag'), {\Phi''}(\ag''))
\\
& = d^*\big((\Phi^*,y),\, ({\Phi'}^*,y')\big)
+d^*\big(({\Phi'}^*,y'),\,({\Phi''}^*,y'')\big)
\end{split}
\end{equation*}
Thus, $d^*$ is a pseudo-metric on $X^*$.

Let us consider the mapping:
\begin{equation*}
\begin{split}
\qquad i:X \to X^* \qquad \mbox{such}& \mbox{ that } x \mapsto
i(x)=(\und x^*, f(x))
\\
& \mbox{where } \, \und x: \{ 0 \} \to X \\
& \mbox{is the constant net defined by }\, \und x(0)=x
\end{split}
\end{equation*}
This mapping is evidently injective. Furthermore, for any $x,x' \in
X$ we have that
\begin{equation*}
\begin{split}
d^*(i(x),i(x')) & = d^*\big((\und x^* ,f(x)),\,( \und x'^*,f(x'))\big)
\\
& = \lim_{(\ag,\ag')\in D\times D'} d(\und x(\ag),\und x'(\ag'))
\\
& = \lim_{(\ag,\ag')\in D\times D'} d(x, x')
\\
& = d(x,x') \, .
\end{split}
\end{equation*}
So, $i:X\to X'$ is an isometric embedding between the pseudo-metric
spaces $(X,d)$ and $(X^*,d^*)$ and we can identify any $x$ with
the corresponding equivalence class $\und x^*$ of the constant net $\und x$.
\\
In particular, $d^*|_{X\times X} \equiv d^*|_{i(X)\times i(X)}=d$.

Let us define a mapping
$$ f^* : X^* \to Y \quad \mbox{ such that } \quad (\Phi^* , y)
\mapsto f^*\big( \Phi^* , y \big) =y.$$
Evidently, for every $x \in X$, we have
$(f^*\circ i)(x)=f^*(i(x))=f^*\big( (\und x^* , \,
f(x)\big)=f(x)$, i.e. that
$$f^*\circ i = f$$
or (identifying $x$ with its isometric image $i(x)$) that
$$f^*|_X = f .$$

Now, for any fixed $y \in Y$ take $(\Phi^*, y), ({\Phi'}^* ,y') \in X^*$ such that
$f^*(\Phi^*, z)=f({\Phi'}^*, z')=y$ and
$d^*\big( (\Phi^*, z) , \, ({\Phi'}^*, z')\big)=0$,
it follows that $z=z'=y$ and $\lim_{(\ag,\ag')\in D\times D'}
d\big( \Phi(\ag), {\Phi'}(\ag')\big)=0$.
Hence, for every $\eg >0$, there exist some $\cg \in D$ and $\cg' \in D'$ such
that for any $\ag>\cg$ and $\ag'>\cg'$, $d\big( \Phi(\ag),
{\Phi'}(\ag')\big)<\eg$. Since $\Phi, {\Phi'} \in C_y$, it
follows that $\Phi \sim_y {\Phi'}$ and so that $(\Phi^*, z)=({\Phi'}^*,z')$.
\\
This proves that, for every $y \in Y$, the restriction
$d^*|_{f^*\s(\{y\})\times f^*\s(\{y\})}$
is a metric and so that $d^*$ is a metric on the mapping
$f^*:X^* \to Y$.

Let us note that the mapping $i:X\to X^*$ is continuous with
respect to the topologies $\tau(f,d)$ and $\tau(f^*,d^*)$.
\\
In fact, for every $N \in \cB(f^*,d^*)$,
there are some $U\in \tau(^*)$ and $V \in \tau (Y)$ such that
$N=U\cap f\s(V)$ and
$$i\s(N)=i\s(U)\cap i\s\big(f^*\s(V)\big) = i\s(U) \cap f\s(V)$$
which is an open set in $\tau(f,d)$ because
the mapping $i$ is evidently continuous with respect to
$\tau(d)$ and $\tau(d^*)$.

To prove that $i(X)$ is dense in $X^*$ with respect to $\tau(f^*,d^*)$
we have to show that for every $W \in \cB (f^*,d^*)\setminus \{ \emptyset \}$ there exists
some $x \in X$ such that $i(x) \in W$.
Let $(\Phi^*, y) \in W=U\cap f^*\s(V)$ with $U \in \tau(d^*)$ and
$V \in \tau (Y)$.
Then, there is some $\eg >0$ such that
$B_{d^*}\big((\Phi^*,y),\eg\big)\sub U$.
Since $\Phi$ is a Cauchy net tied to $y$, there exists some
$\cg\in D$ such that for every $\ag,\bg \> \cg$,
$d(\Phi(\ag),\Phi(\bg))<\frac{\eg}{2}$ and $\Phi(\ag) \in f\s(V)$.
So, took $x=\Phi(\bg)$ for some $\bg>\cg$ we have that
$x\in f\s(V)$ implies $f(x)\in V$, $f^*(i(x))\in V$ and $i(x)\in
f^*\s(V)$.
Since, for every $\ag >\cg$ we have that
$d(\Phi(\ag),x)<\frac{\eg}{2}$ it also follows that
$$d^*\big( (\Phi^*,y),i(x)\big) = \lim_{\ag \in D} d\big(
\Phi(\ag),x\big)<\eg$$
and so that $i(x)\in B_{d^*}\big((\Phi^*,y),\eg\big)$.
This proves that $i(x)\in U \cap f\s(V) \sub W$.

Then $i:X \to X^*$ is a isometric dense embedding from the mapping
$f:X\to Y$ to $f^* : X^* \to Y$.

Finally, we prove that the metric mapping $(f^*, d^*)$ is
complete.
\\
Let $y \in Y$ and $\Psi : E \to X^*$ be a Cauchy net on $X^*$ tied to $y$.
Then, for every $\lm \in E$, $\Psi(\lm)=\big(\Phi^*_\lm ,y_\lm \big)$
where $y_\lm \in Y$ and $\Phi_\lm : D_y \to X$ is a Cauchy net tied at $y_\lm$.

For any fixed $\lm \in E$ and every $n \in \NNs$, since $f\s(O)$ is dense in
$f^*\s(O)$, there exists some $x_{\lm n} \in f\s(O)$ such that
$d^*\big(\Psi(\lm), i(x_{\lm n})\big)<\frac 1n$.

Now, the net $\LM : E \times \NNs \to X$ defined by $\LM(\lm, n)=x_{\lm n}$ (for every
$(\lm , n)$ in the directed set $\LM \times \NNs$)
is evidently tied to $y$ and it is also a Cauchy net.
In fact, for every $\eg >0$, since $\Psi : E \to X^*$ is Cauchy,
there exists some $\xi \in E$ such that for every
$\lm,\mu\ge \xi$,
\begin{equation*}
d^*(\big( \Psi(\lm), \Psi(\mu)\big) < \frac{\eg}{3} .
\end{equation*}
Then, fixed $n_0 \in \NNs$ such that $n_0 > \frac{3}{\eg}$, for every
$(\lm,l),(\mu,m) \in \LM \times \NNs$ such that $(\lm,l),(\mu,m)> (\xi, n_0)$
we have that
\begin{equation*}
\begin{split}
d\big(\LM(\lm, l), \LM(\mu, m) \big)
& = d(x_{\lm l},x_{\mu m})
\\
& = d^*\big( i(x_{\lm l}),i(x_{\mu m})\big)
\\
& \le d^*\big( i(x_{\lm l}),\Psi(\lm )\big)
+ d^*\big( \Psi(\lm), \Psi(\mu)\big)
+ d^*\big( \Psi(\mu),i(x_{\mu m})\big)
\\
& < \frac 1l + \frac{\eg}{3} + \frac 1m
\\
& \le \frac{1}{n_0} + \frac{\eg}{3} + \frac{1}{n_0}
\\
& < \frac{\eg}{3} + \frac{\eg}{3} + \frac{\eg}{3}
\\
& \le \eg \, .
\end{split}
\end{equation*}
So, $\LM \in \cC_y$ and $(\LM^*, y) \in X^*$.

Since, it is clear that $f^*(\LM^*, y)=y$, to finish the proof
it suffices to show that the net $\Psi:E\to
X^*$ converges to $(\LM^*, y)$.
In fact, for every $\eg >0$, since $\LM : E \times \NNs \to X$
is a Cauchy net, there exists some $(\xi, n_0) \in E\times \NNs$
such that for every
$(\lm, l), (\mu, m) \ge (\xi, n_0)$,
$\,\, \ds d(\big( \LM(\lm, l), \, \LM(\mu, m)\big) <
\frac{\eg}{2}$.
\\
Hence, for every $\lm \ge \xi$, we have that
\begin{equation*}
\begin{split}
d^*\big( \Psi(\lm), (\LM^* , y) \big) & \le \eg
\end{split}
\end{equation*}
which concludes our proof.
\qed

{\noindent {\it Key words and phrases:}} complete metric mapping,
completion of a mapping, Cauchy net tied at a point.

\vskip 2pt

{\noindent {\it AMS Subject Classification:}} Primary 54C05, 54C10, 54C20.

\vskip 10pt

\parskip=0pt {\noindent {\sc Giorgio NORDO \newline
MIFT - Dipartimento di Scienze Matematiche e Informatiche, Scienze Fisiche e Scienze della Terra,
Universit\`{a} di Messina\\
Viale F. Stagno D'Alcontres, 31 - Contrada Papardo\\
Sant'Agata -- 98166 Messina, Italy.\\
E-mail: {\tt giorgio.nordo@unime.it}

\end{document}